\documentclass[12pt,twoside]{article} 
\usepackage[latin1]{inputenc}
\usepackage[T1]{fontenc}
\usepackage{graphicx}
\usepackage{enumitem}
\usepackage{latexsym}
\usepackage{amssymb}
\usepackage{epsfig}
\usepackage{xcolor}
\usepackage{amsmath}
\usepackage{natbib}
\usepackage{float}
\allowdisplaybreaks 

\numberwithin{equation}{section}

\newcommand{\Var}{\mbox{Var}}
\def\eps{\varepsilon}
\def \R{\mathbb{R}}

\def \Cov{\mbox{Cov}}

\newcommand{\nt}{\lfloor nt\rfloor}

\def\er{\mathbb{R}}

\def\cn{\mathcal{N}}

\def\e{\varepsilon}

\def\beq{\begin{eqnarray*}}
\def\eeq{\end{eqnarray*}}

\usepackage{titling}

\begin{document}

\title{\bf Testing for changes in the error distribution in  functional linear models
}

\author{\textsc{ Natalie Neumeyer}\thanks{corresponding author; e-mail: natalie.neumeyer@uni-hamburg.de; orcid 0000-0002-6649-1135} \qquad \textsc{Leonie Selk}\thanks{e-mail: leonie.selk@uni-hamburg.de; orcid 0000-0003-0459-4953}\\  Universit\"at Hamburg, Fachbereich Mathematik}
\thanksmarkseries{arabic}

\maketitle

\newtheorem{theo}{Theorem}[section]
\newtheorem{lemma}[theo]{Lemma}
\newtheorem{cor}[theo]{Corollary}
\newtheorem{rem}[theo]{Remark}
\newtheorem{prop}[theo]{Proposition}
\newtheorem{defin}[theo]{Definition}
\newtheorem{example}[theo]{Example}
\newtheorem{Assumption}{Assumption}

\begin{abstract}
We consider linear models with scalar responses and covariates from a separable Hilbert space. The aim is to detect change points in the error distribution, based on sequential residual empirical distribution functions. Expansions for those estimated functions are  more challenging in models with infinite-dimensional covariates than in regression models with scalar or vector-valued covariates  due to a slower rate of convergence of the parameter estimators. Yet the suggested change point test is asymptotically distribution-free and consistent for one-change point alternatives. In the latter case we also show consistency of a change point estimator. 
\end{abstract}

MSC 2020 Classification: Primary  62R10
Secondary 62G10, 62G30

Keywords and Phrases: change-points, functional data analysis, regularized function estimators, regression, residual processes


\section{Introduction}

We consider a functional linear model $Y=\alpha+\langle X,\beta\rangle+\eps$ with scalar response $Y$ and covariates $X$ from a separable Hilbert space, e.g.\ $L^2([0,1])$. 
Structural changes in the distribution can appear, even when the parameters $\alpha$ and $\beta$ do not change. For this reason we focus on detecting changes in the error distribution.
If the errors were observable one could use the classical test (and change point estimators) based on the difference of the sequential empirical distribution functions of the first $\lfloor nt\rfloor$ and the last $n-\lfloor nt\rfloor$ error terms from a sample of $n$ observations, see \cite{CsorgoEtal1997}, 
\cite{Picard1985}, \cite{Carlstein1988}, \cite{Dumbgen1991}, \cite{HarizEtal2005} and \cite{HarizEtal2007}. In a regression model those tests have to be 
 based on estimated residuals $\hat \eps=Y-\hat\alpha-\langle X,\hat\beta\rangle$.
 Similar tests have been considered by \cite{Bai1994} in the context of ARMA-models, by \cite{Koul1996} in the context of nonlinear time series, by \cite{Ling1998} for nonstationary autoregressive models, and by 
\cite{NeumeyerKeilegom2009} and
\cite{SelkNeumeyer2013} for nonparametric independent and time series regression models. Typically the asymptotic distribution is derived using asymptotic expansions of residual-based empirical distribution functions. For models with functional covariates those expansions can be  problematic because inner products $\langle X,\hat\beta-\beta\rangle$ appear and those can have a slow rate of convergence (see \cite{CardotEtal2007}, \cite{ShangCheng2015}, and \cite{YeonEtal2023}). However, we show that under very simple non-restrictive assumptions those terms cancel for the suggested change point test statistic and thus the asymptotic distribution is the same as based on true (unobserved) errors. 

Change point testing and estimation for functional data, and for the parameter in functional linear models have been considered in the literature, but not for the error distribution. Tests for changes in the functional mean and in the parameter function of autoregressive models are considered in chapters 6 and 14 in \cite{HorvathKokoszka}. 
In \cite{BerkesEtal2009} a CUSUM testing procedure to detect a change in the mean of functional observations is proposed. They apply projections on principal components of the data to estimate the mean.
\cite{AueEtal2009} extend this result and introduce an estimator for the change point in this model and derive its limit distribution. \cite{AstonKirch2012} consider the same type of model with epidemic changes and dependent data.
In \cite{AueEtal2018} structural breaks in the mean of functional observations are detected and dated without the application of dimension reduction techniques (as functional principal component analysis). 
\cite{AueEtal2014} propose a monitoring procedure to detect structural changes in functional linear models with functional response, allowing for dependence in the data, including functional autoregressive processes. They test for a change in the regression operator, which is the analogue to our  $\beta$, based on functional principal component analysis.
A linear regression model with scalar response is considered in \cite{HorvathEtal2023} who propose a tests for the detection of multiple change points in the regression parameter. The regressors in their model can be functional and can include lagged values of the response.



The paper is organized as follows. In Section 2 we define the test statistic and present model assumptions to obtain the asymptotically distribution-free hypothesis test. In Section 3 we discuss the assumptions on the parameter estimators and some examples. In Section 4 consistency of the test as well as of a change point estimator is considered in the context of one change point. Finite sample properties are shown in Section 5. Section 6 concludes the paper, in particular with an outlook on goodness-of-fit testing. The proofs are given in the appendix.

\section{Model, test statistic and main result under the null}\label{sec-model}

Let $\mathcal H$ be a separable Hilbert space with inner product $\langle \cdot,\cdot\rangle$, corresponding norm $\|\cdot\|$ and Borel-sigma field.  Let $(X_i,Y_i)$, $i=1,\dots,n$, be an independent sample of ($\mathcal{H}\times \mathbb{R}$)-valued random variables defined on the same probability space with probability measure $\mathbb{P}$. The data  are modeled as functional linear model
\begin{equation*}
Y_i=\alpha+ \langle X_i,\beta\rangle+ \eps_i,\quad i=1,\dots,n,
\end{equation*}
 with scalar response $Y_i$ and $\mathcal H$-valued covariate $X_i$, and with  parameters $\alpha\in\R$, $\beta\in H$. The covariates $X_1,\dots,X_n$ are assumed to be iid with $E \Vert X_i \Vert< \infty$, and the errors $\eps_1,\dots,\eps_n$ are independent, centered, and independent of the covariates. 
 Our aim is to test for change-points in the error distribution. In this section we consider the test statistic under the null hypothesis, where the errors are identically distributed. 

Let $\hat\alpha$ and $\hat\beta$ denote estimators for the parameters $\alpha\in\mathbb{R}$ and $\beta\in \mathcal H$. We build residuals 
$\hat\eps_i=Y_i-\hat\alpha-\langle X_i,\hat\beta\rangle$, $i=1\dots,n$. The test statistic
$$T_n=\sup_{t\in[0,1]}\sup_{z\in\R}|\hat G_n(t,z)|$$
based on the process
$$\hat G_n(t,z)=\frac{\nt \left(n-\nt\right)}{n^{3/2}}\left(\hat F_{\nt}(z)-\tilde F_{\nt}(z)\right),$$
compares for each $k=1,\dots,n-1$ the   empirical distribution functions
$$\hat F_k(z)=\frac{1}{k}\sum_{i=1}^k I\{\hat\eps_i\leq z\},\qquad \tilde F_k(z)=\frac{1}{n-k}\sum_{i=k+1}^n I\{\hat\eps_i\leq z\}$$
of the first $k$ and last $n-k$ residuals, respectively. 
Note that one can write 
$$\hat G_n(t,z)=\frac{\nt}{n^{1/2}}\left(\hat F_{\nt}(z)- \hat F_{n}(z)\right).$$

For the asymptotic distribution of the test statistic under the null hypothesis we assume the following conditions. Let
$P$ denotes the distribution of $(X_1,\eps_1)$. 

\begin{itemize}
\item[\bf (a.1)] $|\hat\alpha-\alpha|=o_{\mathbb{P}}(1)$, $\Vert\hat\beta-\beta\Vert=o_{\mathbb{P}}
(1)$
\item[\bf (a.2)] Let $\eps_1,\dots,\eps_n$ be independent and identically distributed with cdf $F$ that is H\"{o}lder-continuous of order $\gamma\in (0,1]$ with H\"{o}lder-constant $c$.
\item[\bf (a.3)]$\mathbb{P}\big(\hat \beta-\beta\in\mathcal{B}\big)\to 1$ as $n\to\infty$ for 
a  class $\mathcal{B}\subset\mathcal H$ such that the function class
$$\mathcal{F}= \left\{(x,e)\mapsto I\{e\leq v+\langle x,b\rangle\}\mid v\in\R, b\in \mathcal{B}\right\}$$
is $P$-Donsker. 
\end{itemize}

\begin{rem}\label{RemarkAssumptions}\rm
The assumptions are very mild and in particular less restrictive than typical assumptions for asymptotic distribution of residual-based empirical processes, even for finite-dimensional covariates. In assumption (a.1)  only consistency is needed, no rates of convergence. Typically in the literature about residual-based procedures a bounded error density is assumed, see e.g.\ \cite{AkritasKeilegom2001}. Then (a.2) is fulfilled for $\gamma=1$, but (a.2) is less restrictive in the cases $\gamma\in (0,1)$. Suitable conditions for the general assumption (a.3) are discussed in Section 3. One possibility for $\mathcal{H}=L^2([0,1])$ is to assume smoothness of $\beta$ which is a typical assumption. If $\gamma\in (\frac12,1]$ in assumption (a.2), and $\beta$ is in a Sobolev-space with third derivatives, (a.2) holds for the estimator $\hat\beta$ from \cite{YuanCai2010}. This estimator can also be applied for smaller $\gamma$ in (a.2) if higher smoothness of $\beta$ is assumed.   
\end{rem}

Define the process $G_n$ as $\hat G_n$, but based on the true errors  instead of residuals, i.e.
$$G_n(t,z)=\frac{\nt}{n^{1/2}}\left( F_{\nt}(z)- F_{n}(z)\right)$$
with 
\begin{equation}\label{F_nt}
F_{\nt}(z)=\frac{1}{\nt}\sum_{i=1}^{\nt} I\{\eps_i\leq z\}.
\end{equation}
Further let $G$ be a completely tucked Brownian sheet, i.e.\ a centered Gaussian process on $[0,1]^2$ with covariance structure
$$\Cov( G(s,u),G(t,v))=(s\wedge t-st)(u\wedge v-uv).$$

\begin{theo}\label{theo1}
Under the assumptions (a.1)--(a.3), 
\begin{equation}\label{hat-G_n-G_n}
\sup_{t\in [0,1]}\sup_{z\in\R}|\hat G_n(t,z)-G_n(t,z)|=o_{\mathbb{P}}(1),
\end{equation}
and thus the process $(\hat G_n(t,z))_{t\in[0,1],z\in\R}$ converges weakly to $( G(t,F(z)))_{t\in[0,1],z\in\R}$. 
\end{theo}

\noindent The proof of (\ref{hat-G_n-G_n}) in the theorem is given in the appendix.  The weak convergence of $G_n$ is a classical result, see \cite{BickelWichura1971} and \cite{ShorackWellner}. 
With the continuous mapping theorem one obtains the asymptotic distribution of the test statistic $T_n$ under the null hypothesis of no change-point, which is the distribution of $T=\sup_{t,u\in [0,1]}|G(t,u)|$ because $F$ is continuous. 
 The test statistic is asymptotically distribution-free with the same limit distribution as for corresponding changepoint tests based on iid observations (not residuals).
Let $\bar\alpha\in (0,1)$ and $q$ be the $(1-\bar\alpha)$-quantile of $T$. Then the test that rejects the null hypothesis if $T_n>q$ has asymptotic level $\bar\alpha$.
Consistency is considered in Section 4.   

\begin{rem}\label{RemarkCvM}\rm
The choice of $T_n$ as a Kolmogorov-Smirnov type test statistic is not mandatory.
In principle, any continuous functional of the process $\hat G_n$ can be considered. The most common ones, besides $T_{n}$,
are of Cram\'{e}r-von-Mises type, e.\,g.\ $T_{n,2}= \sup_{t\in[0,1]}\int_\mathbb{R}|\hat G_n(t,z)|^2dF(z)$ or $T_{n,3}=\int_0^1\int_{\mathbb{R}} |\hat G_n(t,z)|^2dF(z)dt$. The asymptotic distribution of these test statistics under the null hypothesis also follows from Theorem \ref{theo1} and with
$$\int_\mathbb{R} |\hat G_n(t,z)|^2dF(z)\to \int_\mathbb{R}|G(t,F(z))|^2dF(z)=\int_0^1|G(t,x)|^2dx,$$ 
and thus these test statistics are asymptotically distribution-free as well.  However,  $T_{n,2}$ and $T_{n,3}$ contain the unknown quantity $F$ and must therefore be modified in order to be applied. This can be done by replacing the integral with the sample mean:  $\tilde T_{n,2}=\sup_{t\in[0,1]}\frac 1n\sum_{i=1}^n |\hat G_n(t,\hat\varepsilon_i)|^2$ and  $\tilde T_{n,3}=\int_0^1\frac 1n\sum_{i=1}^n |\hat G_n(t,\hat\varepsilon_i)|^2dt$.
\end{rem}

\section{Discussion of assumptions and examples}

To show validity of the Donsker-class assumption (a.3) there are sufficient conditions on covering numbers or bracketing numbers. We discuss some specific conditions on the class $\mathcal{B}$, examples for Hilbert spaces $\mathcal{H}$, and estimators for the parameter function $\beta$ that fulfill the conditions.

\subsection{VC-class condition}

Assumption (a.3) can be derived from a VC-function class condition  formulated as follows.
Assume that
$\mathbb{P}\big(\hat \beta-\beta\in\mathcal{B}\big)\to 1$ as $n\to\infty$ for 
a   class $\mathcal{B}\subset\mathcal H$ such that the class of maps
\begin{equation}\label{VC}
\{\mathcal{H}\to\mathbb{R}, x \mapsto \langle x,b\rangle+v\mid b\in\mathcal{B},v\in\mathbb{R}\} 
\end{equation}
 is  a VC-subgraph class. By definition then 
$$\left\{ \{(x,e)\in\mathcal{H}\times\mathbb{R}\mid e\leq \langle x,b\rangle+v\}\mid  b\in\mathcal{B},v\in\mathbb{R} \right\}
$$
is a VC-class of sets. The class  $\mathcal F$ from (a.3) is the class of the corresponding indicator functions
 and (a.3) is fulfilled by Theorems 8.19 and 9.2 in \cite{Kosorok}.

\begin{example} \rm
We consider the Hilbert space $\mathcal{H}=L^2([0,1])$ with inner product $\langle g,h\rangle=\int_0^1 g(t)h(t)\,dt$ and norm $\|g\|=(\int_0^1 g^2(t)\,dt)^{1/2}$. For the parameter function $\beta$ we assume sparsity as in \cite{LeePark2012}. Let $(\phi_j)_{j\in\mathbb{N}}$ be a basis of $\mathcal{H}$ and assume $\beta=\sum_{j\in J}\beta_j\phi_j$ for some finite, but unknown index set $J$. \cite{LeePark2012} consider the estimator $\hat\beta=\sum_{j=1}^k\hat\beta_j\phi_j$ with 
$$(\hat\beta_1,\dots,\hat\beta_k)=\arg\min_{b_1,\dots,b_k\in\mathbb{R}}\left(\frac1n\sum_{i=1}^n \Big(Y_i-\overline{Y}_n-\sum_{j=1}^kb_j\langle X_i-\overline{X}_n,\phi_j\rangle\Big)^2+\sum_{j=1}^k\hat w_j|b_j|\right)^2,$$
where $k$ is a chosen dimension-cut-off, $\hat w_j$ are suitable weights based on initial estimators, and $\overline{Y}_n=\frac1n\sum_{i=1}^n Y_i$, $\overline{X}_n=\frac1n\sum_{i=1}^n X_i$. Further, $\hat\alpha=\overline{Y}_n-\langle\hat\beta,\overline{X}_n\rangle$.
Under suitable assumptions, in particular $E\|X\|^2<\infty$, and $k$ is larger than the largest index in $J$, \cite{LeePark2012} show in their Theorem 2 that $\mathbb{P}(\hat\beta_j=0\mbox{ for }j\not\in J)\to 1$ for $n\to\infty$. Thus we can set
$$\mathcal{B}=\left\{\sum_{j\in J}b_j\phi_j\;\Big|\; b_j\in\mathbb{R}\,\forall j\in J\right\}$$
and obtain $\mathbb{P}(\hat\beta-\beta\in\mathcal{B})\to 1$ for $n\to\infty$. Further, 
 the class of maps in (\ref{VC}), i.e.\
$$\left\{\mathcal{H}\to\mathbb{R}, x \mapsto \sum_{j\in J}b_j\langle x,\phi_j\rangle+v\;\Big|\; b_j\in\mathbb{R}\,\forall j\in J, v\in\mathbb{R}\right\} $$ 
 is a finite dimensional vector space and thus a VC-class, see Lemma 2.6.15 in \cite{VaartWellner}. Then as discussed above validity of (a.3) follows. Furthermore, from Theorem 2 in \cite{LeePark2012} it also follows that our assumption (a.1) is fulfilled, and thus under assumption (a.2) the assertion of Theorem \ref{theo1} holds. 
\end{example}

\subsection{Bracketing number condition}

In this subsection we assume that $\mathcal H$ is a separable Hilbert space of real-valued functions (or vectors with real components) and the inner product is increasing in the sense that from $h\leq g$ (pointwise for functions; componentwise for vectors) it follows that $\langle h,x\rangle \leq \langle g,x\rangle$ for all $x\in\mathcal H$ with $x\geq 0$. Then assumption (a.3) can be replaced by the condition in the next lemma.

\begin{lemma}\label{lem1}
Assume (a.1), (a.2) and $\mathbb{P}\big(\hat \beta-\beta\in\mathcal{B}\big)\to 1$ as $n\to\infty$ for 
a  function class $\mathcal{B}\subset\mathcal H$ such that the bracketing number fulfills
$ \log N_{[\,]}(\mathcal{B},\epsilon,\|\cdot\|)\leq K/\epsilon^{1/k}$ 
for some $k>1/\gamma$. Here $\gamma$ is the H\"{o}lder-order from assumption (a.2). 
Then assumption (a.3) holds.  
\end{lemma}

The proof is given in the appendix.

\begin{example} \rm\label{ex-bracketing}
We consider the Hilbert space $\mathcal{H}=L^2([0,1])$ with inner product $\langle g,h\rangle=\int_0^1 g(t)h(t)\,dt$ and norm $\|g\|=(\int_0^1 g^2(t)\,dt)^{1/2}$. 
We assume $\beta\in \mathcal{W}_2^m([0,1])$ for some $m>2$ and the Sobolev-space
\begin{eqnarray*}
\mathcal{W}_2^m([0,1])& &= \big\{b:[0,1]\to\mathbb{R}\mid 
b^{(j)}\mbox{ is absolutely continuous for }j=0,\dots,m-1,\\
&&\qquad\qquad\qquad\qquad\mbox{ and }\|b^{(m)}\|<\infty
\big\},
\end{eqnarray*}
where $b^{(0)}=b$, and $b^{(j)}$ denotes the $j$-th derivative of $b$, $j\geq 1$.
We consider the regularized estimators in  \cite{YuanCai2010}, i.e.
\[
\big(\hat\alpha,\hat\beta \big)=\arg\min_{a\in\mathbb{R},b\in \mathcal{W}_2^m([0,1])} \left\{\frac1n\sum_{i=1}^n 
\Big(Y_{i}- \big(a+\langle X_i,b\rangle \big)\Big)^2+\lambda_n\big\|b^{(m)}\big\|^2\right\}
\]
for a suitable positive sequence $\lambda_n$ converging to zero.
 Convergence rates of $\hat\beta$ and its derivatives can be found in Corollaries 10 and 11 in \cite{YuanCai2010}. Under suitable assumptions one obtains $\big\|\hat \beta^{(j)}-\beta^{(j)}\big\|=o_{\mathbb{P}}(1)$ for $j=0,1,2$, and thus $\mathbb{P}(\hat\beta-\beta\in\mathcal{B})\to 1$ for the function class
\[
\mathcal{B}=\big\{b\in \mathcal{W}_2^{2}\big([0,1]\big): \|b\| + \|b^{(2)}\|\leq 1\big\}. 
\] 
By Corollary 4.3.38 
in \cite{GineNickl} and Lemma 9.21 in \cite{Kosorok}  the class $\mathcal B$ fulfills the bracketing number condition in Lemma \ref{lem1} for $k=2$.
 Thus the assumptions (a.1)--(a.3) are fulfilled if $F$ is H\"{o}lder-continuous of order $\gamma\in (\frac12,1]$. 
 Less restrictive assumptions on $F$, i.e.\ $\gamma\leq \frac12$, require for this concept higher smoothness of $\beta$.

\end{example}

\section{Fixed one-change point alternative: consistency of the test and change point estimator}

 In this section we consider fixed alternatives with one change point at index $k_n^*=\lfloor n\vartheta^*\rfloor$ with $\vartheta^*\in (0,1)$. We write the functional linear model as in Section 2 under the following assumption. 
 
\begin{itemize}
\item[\bf (a.2)'] Assume $\eps_{1},\dots,\eps_{k^*_n}$ are iid with cdf $F_1$, and $\eps_{k_n^*+1},\dots,\eps_{n}$ are iid with cdf $F_2\neq F_1$. Let  $F_1$ and $F_2$ be H\"{o}lder-continuous of order $\gamma_1,\gamma_2\in (0,1]$ with H\"{o}lder-constant $c_1,c_2$, respectively.
\end{itemize}

Let further $P_1$ denote the distribution of $(X_1,\eps_1)$ (before the change) and $P_2$ denote the distribution of $(X_n,\eps_n)$ (after the change).
For the empirical distribution functions $\hat F_k$ and $\tilde F_k$ as in Section 2 we obtain the following asymptotic result.

\begin{lemma}\label{lem-consistency}
Under assumptions (a.1) and (a.2)' and if (a.3) is valid for $P=P_1$ and $P=P_2$, it holds that
$$\sup_{z\in\mathbb{R}}|\hat F_{k_n^*}(z)-F_1(z)|=o_{\mathbb{P}}(1)\mbox{ and } \sup_{z\in\mathbb{R}}|\tilde F_{k_n^*}(z)-F_2(z)|=o_{\mathbb{P}}(1).$$
\end{lemma}
The proof is given in the appendix. 
Now note that 
$$\frac{T_n}{n^{1/2}}\geq \frac{k_n^*(n-k_n^*)}{n^{2}}\sup_{z\in\mathbb{R}}\left| \hat F_{k_n^*}(z)-\tilde F_{k_n^*}(z)\right|,$$
and by Lemma \ref{lem-consistency} the right hand side converges in probability to the positive constant 
$$\vartheta^*(1-\vartheta^*)\sup_{z\in\mathbb{R}}\left| F_1(z)- F_2(z)\right|.$$
From this it follows that tests that reject the null hypothesis of no change-point if $T_n>q$ for some $q>0$ (see Section 2) are consistent.



The estimator for the change point $\vartheta^*$ is based on the process $\hat G_n$ and is defined as
$$\hat\vartheta_n=\min\Big\{t:\sup_{z\in\mathbb{R}}|\hat G_n(t,z)|=\sup_{t'\in[0,1]}\sup_{z\in\mathbb{R}}|\hat G_n(t',z)|\Big\}.$$

\begin{lemma}\label{lem-est-consistency}
Under assumptions (a.1), (a.2)' and if (a.3) holds for $P=P_1$ and $P=P_2$, the change point estimator is consistent, i.\,e.
$$|\hat\vartheta_n-\vartheta^*|=o_{\mathbb{P}}(1).$$
\end{lemma}
The proof is given in the appendix.

\section{Finite sample properties}

We consider the Hilbert space $\mathcal{H}=L^2([0,1])$. 
For $i=1,\ldots,n$ the functional observations $X_{i}(t)$, $t \in [0,1]$, are generated according to
$$X_{i}(t)=\frac 12\sum_{l=1}^5\Big(B_{i,l}\sin\left( t(5-B_{i,l})2\pi\right)-M_{i,l}-E[B_{i,l}\sin\left( (5-B_{i,l})2\pi\right)-M_{i,l}]\Big),$$ 
where $B_{i,l}\sim\mathcal{U}[0,5]$ and $M_{i,l}\sim\mathcal{U}[0,2\pi]$ for $l=1,\ldots,5$, $i=1,\ldots,n$.  $\mathcal{U}$ stands for the (continuous) uniform distribution. The  functional linear model is built as
$$Y_i=\int X_{i}(t)\gamma_{3,\frac 13}(t)dt+\varepsilon_i,$$
where the coefficient function $\gamma_{a,b}(t)=b^a/\Gamma(a)t^{a-1}e^{-bt}I\{t>0\}$ is the density of the Gamma distribution. Furthermore, we assume that each $X_{i}$ is observed on a dense, equidistant grid of 300 evaluation points.

The parameter estimators are the regularized estimators described in Example \ref{ex-bracketing} with $m=3$ and a data-driven tuning parameter $\lambda_n$ chosen by generalized cross-validation as described in \cite{YuanCai2010}.

We model three similar types of change points, such that
\[ \e_1,\ldots,\e_{\lfloor \frac n2\rfloor}\sim \cn(0,1),\ \e_{\lfloor \frac n2\rfloor+1},\ldots,\e_n\sim\tilde{F}_{1,\delta}\ (\text{respectively}\ \tilde{F}_{2,\delta},\ \tilde{F}_{3,\delta}),\]
where $\tilde{F}_{1,\delta}$, $\tilde{F}_{2,\delta}$, $\tilde{F}_{3,\delta}$ have in common that the mean remains zero and the variance remains one. In particular
\begin{itemize}
\item $\tilde{F}_{1,\delta}$ is the distribution function of a random variable that is $\cn(-2\delta,1)$ distributed with probability $0.5$ and $\cn(2\delta,1)$ distributed with probability $0.5$. 
\item $\tilde{F}_{2,\delta}$  is the distribution function of a random variable that is $\cn(0,(1-\delta)^2)$ distributed with probability $0.5$ and $\cn(0,2-(1-\delta)^2)$ distributed with probability $0.5$.
 
\item  $\tilde{F}_{3,\delta}$ is the ``skew-normal''-distribution $$SN\left(-\sqrt{\frac{2\pi\left(\left(10\delta\right)^2+\left(10\delta\right)^4\right)}{\pi^2+\left(2\pi^2-2\pi\right)\cdot\left(10\delta\right)^2+\left(\pi^2-2\pi\right)\cdot\left(10\delta\right)^4}} , \sqrt{\frac{\pi\left(1+\left(10\delta\right)^2\right)}{\pi+(\pi-2)\cdot\left(10\delta\right)^2}},10\delta \right).$$
 A random variable $Z$ is distributed $SN(\lambda_1,\lambda_2,\lambda_3)$ if $Z=\lambda_1 +\lambda_2 \cdot Z_0$ and $Z_0$ has the density $2\phi(x) \Phi(\lambda_3 x)$, where $\phi$ is the density and $\Phi$ is the distribution function of the standard normal distribution (see \cite{AzzaliniCapitanio1999}). 
 The expected value of such a random variable $Z$ is calculated as $\lambda_1+\lambda_2\sqrt{\frac 2\pi}\cdot\frac{\lambda_3}{\sqrt{1+{\lambda_3}^2}}$ and the variance as ${\lambda_2}^2\left(1- \frac2\pi\cdot\frac{{\lambda_3}^2}{1+{\lambda_3}^2}\right)$. This results in the parameters for the distribution after the change point, such that the the expected value of the errors remains 0 and the variance remains 1. 
\end{itemize}
So $\delta=0$ represents the null hypothesis of no change point, and the difference between the distribution before and after the change point grows with $\delta$ in all three cases.

\begin{figure}
\includegraphics[width=0.32\textwidth]{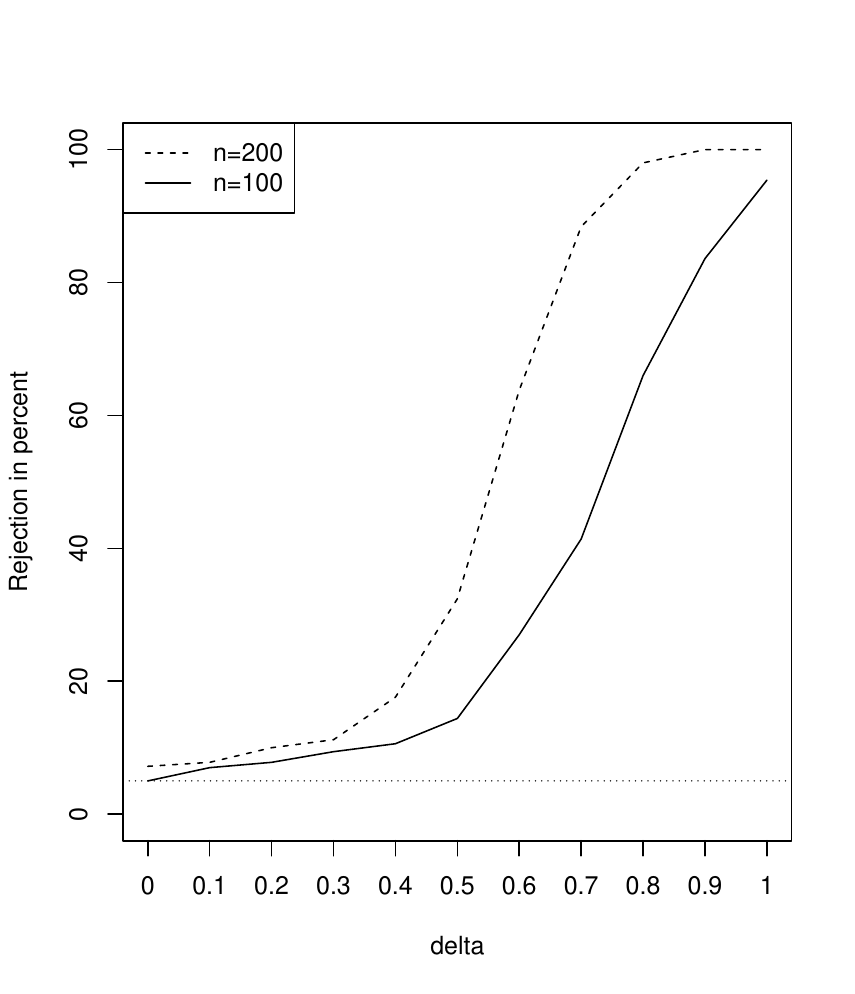} \includegraphics[width=0.32\textwidth]{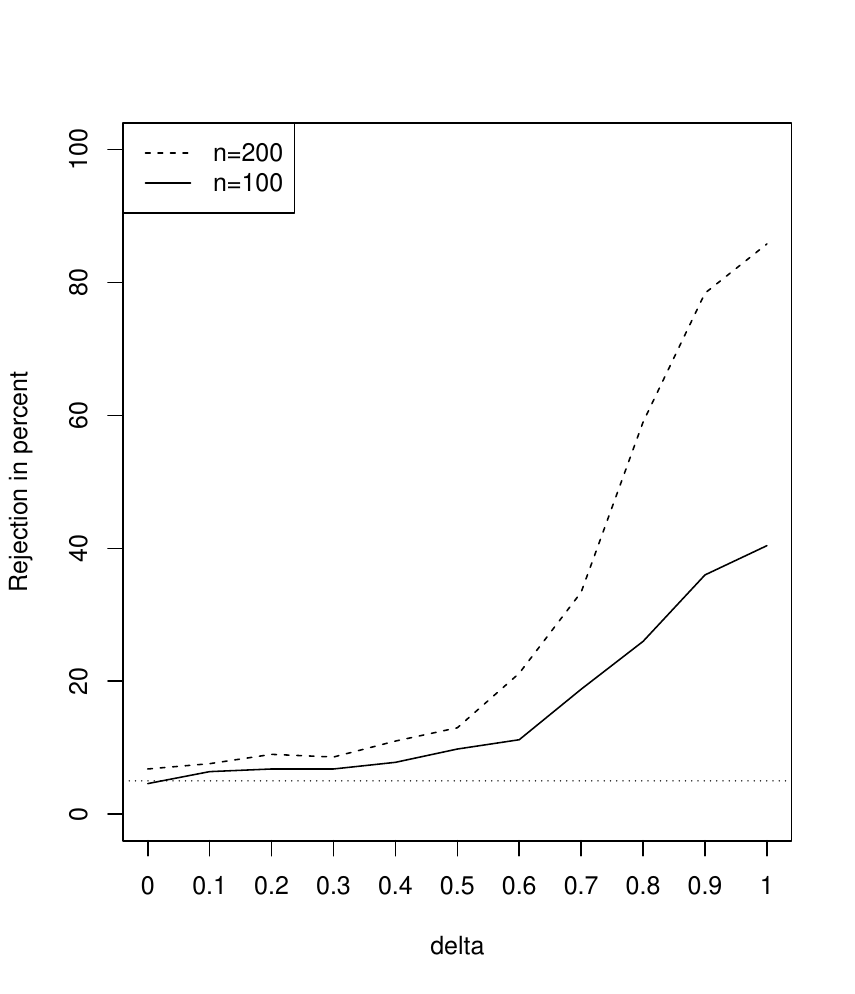} \includegraphics[width=0.32\textwidth]{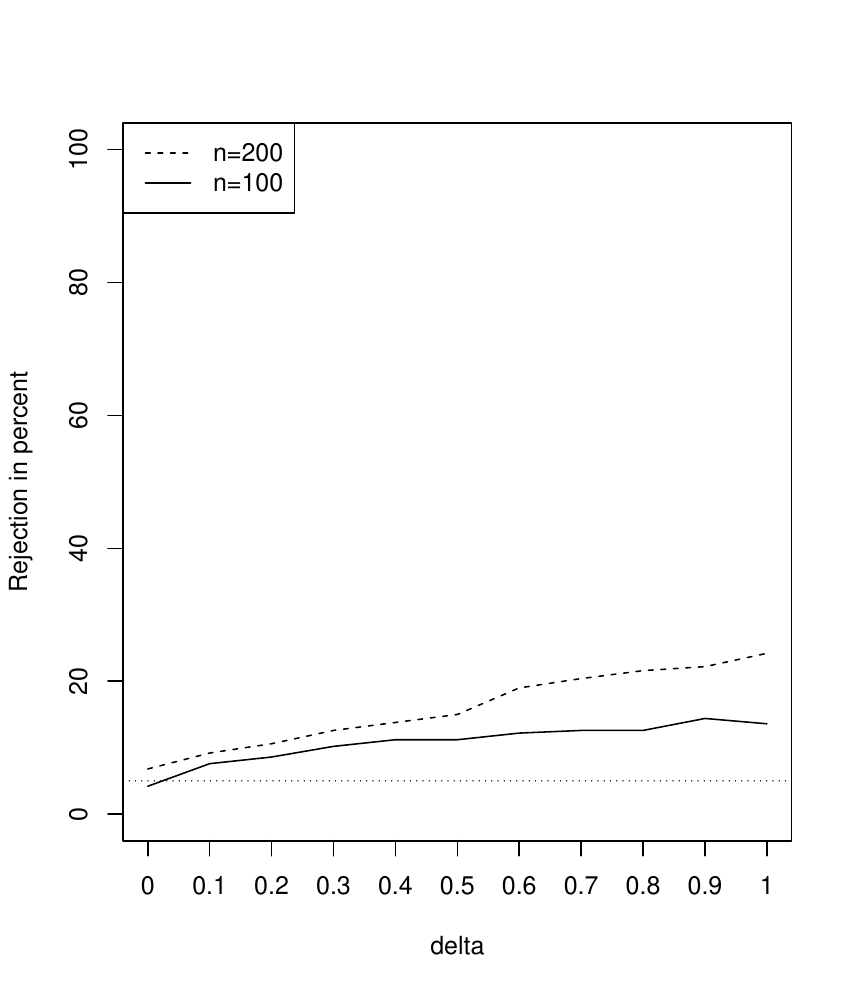}
\caption{Rejection probabilities with a change in the error distribution from $\cn(0,1)$ to $\tilde{F}_{1,\delta}$ (left), to $\tilde{F}_{2,\delta}$ (middle) and to $\tilde{F}_{3,\delta}$ (right). The dotted line marks the $5\%$ level.}
\label{fig:rejprob}
\end{figure}

In Figure \ref{fig:rejprob} the rejection probabilities for 500 repetitions, level 5\% (critical value tabled in  \cite{Picard1985}) and sample sizes $n\in\{100,200\}$ are shown. 
In all three cases it can be seen that the level is approximated well and the power increases for increasing parameter $\delta$ as well as for increasing sample size $n$. In the case of a change in skewness, the increase with $\delta$ is not as pronounced as in the other two cases. This is because the distributions for different values of $\delta$ become more similar as $\delta$ increases.
The same types of changes (from $\cn(0,1)$ to $\tilde{F}_{1,\delta}$ and to $\tilde{F}_{2,\delta}$) were also simulated in \cite{SelkNeumeyer2013} for a real-valued nonparametric autoregression model with lag 1. The results are comparable with an even higher power in the paper at hand.

\begin{figure}
\begin{center}
\includegraphics[width=0.32\textwidth]{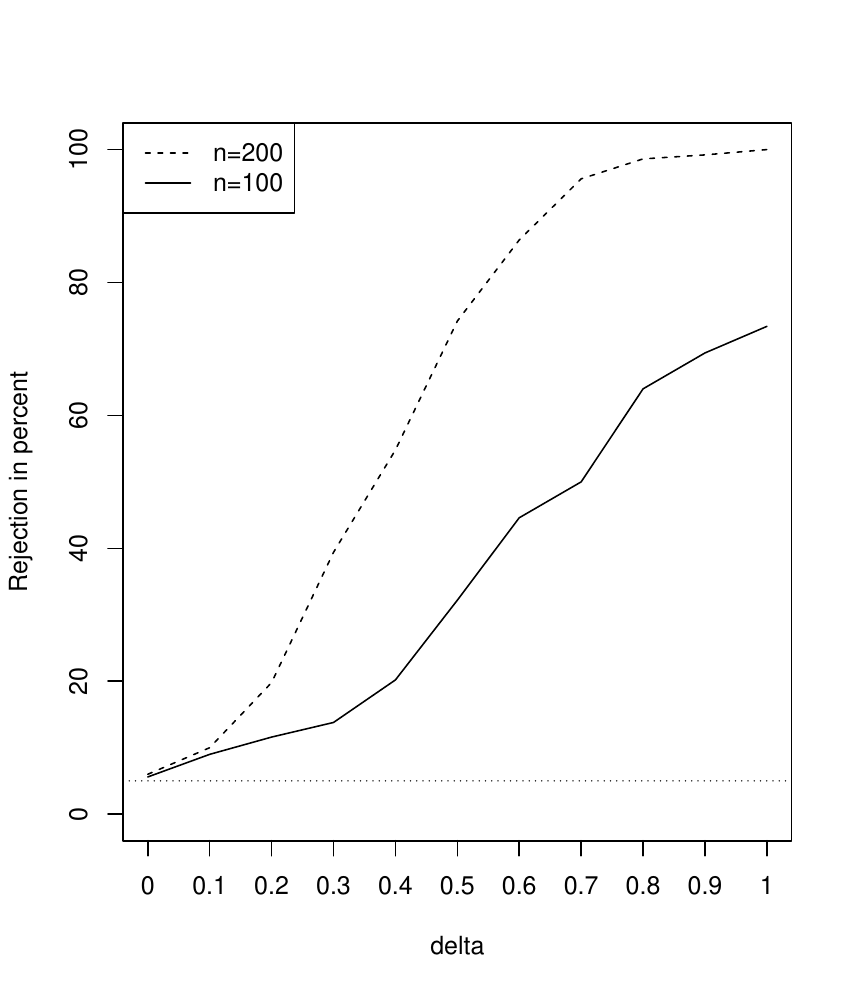}
\end{center}
\caption{Rejection probabilities with a change in the error distribution from $\cn(0,0.5^2)$ to $\cn(0,(0.5+\delta)^2)$. The dotted line marks the $5\%$ level.}
\label{fig:rejprob2}
\end{figure}

In addition, we model a more distinct change, that is
\[ \e_1,\ldots,\e_{\lfloor \frac n2\rfloor}\sim \cn(0,0.5^2),\ \e_{\lfloor \frac n2\rfloor+1},\ldots,\e_n\sim\cn(0,(0.5+\delta)^2).\]
As expected for a change in the variance the power grows faster with increasing $\delta$ than in the other three cases, especially for small $\delta$. The results are shown in Figure \ref{fig:rejprob2}.
This kind of change point was also simulated in \cite{NeumeyerKeilegom2009} for a nonparametric regression model with one-dimensional regressor. The results are similar.

Next to the Kolmogorov-Smirnov type test with test statistic $T_n$, we also applied Cram\'{e}r-von-Mises type tests with test statistics $T_{n,2}$ and $T_{n,3}$. The results are very similar and are not presented here for the sake of brevity.

\section{Concluding remarks}\label{concluding remarks}

To detect structural changes in functional linear models, we considered the classical test by \cite{BickelWichura1971} for a change in the distribution, but based on estimated errors. We gave simple assumptions under which the asymptotic distribution of the test statistic under the null is the same as for iid data. The test as well as the corresponding change point estimators are consistent in one-change point models. The same test can be considered in more complex regression models with functional covariates, e.g.\ a quadratic model as in \cite{BoenteParada2023} or nonparametric models, see \cite{FerratyVieu}. 
We only consider independent data, but testing for change points in the innovation distribution in  times series models that include functional covariates is a very interesting topic. However, the proofs for asymptotic distributions will be  more complicated. In future work we are planning to consider a time series model
		$Y_t=m(X_t)+\varepsilon_t,$
		where $Y_t$ and $\varepsilon_t$ are real-valued and $X_t$ contains a functional part, but can also contain past values $Y_{t-1},\dots,Y_{t-p}$. We presume that the proofs as in  \cite{SelkNeumeyer2013} (for nonparametric autoregression time series with independent errors) and of section 4.2 in  \cite{NeumeyerOmelka2024} (for linear models with finite-dimensional covariates and beta-mixing errors) can be combined with the proofs in the paper at hand to consider change-point tests for the innovation distribution under the assumption that 
		 $(X_t,\varepsilon_t)$, $t\in\mathbb{Z}$, is a strictly stationary beta-mixing time series.

In the proof of Theorem \ref{theo1} we derive an expansion for the sequential residual-based empirical distribution function, 
$$\hat F_{\nt}(z)=F_{\nt}(z)+R_n(z)+o_{\mathbb{P}}(n^{-1/2})$$ 
uniformly in $t\in[0,1],z\in\R$, where
$F_{\nt}$ is defined in (\ref{F_nt}), 
 and the term 
$$R_n(z)=E_X[F(z+\hat\alpha-\alpha+\langle X,\hat\beta-\beta\rangle)]-F(z)$$
 appears from estimating the parameters (see (\ref{Entwicklung-hatF_nt}) in the appendix). 
Here $E_X$ denotes the expectation with respect to $X$, which has the same distribution as $X_i$, but is independent of $\hat\alpha$, $\hat\beta$.  
 For change-point testing  the remainder term $R_n$ cancels when considering the test statistic $T_n$. For other testing procedures, e.\,g.\ goodness-of-fit tests for the error distribution, this typically nonnegligible term is of relevance, see \cite{Koul} and \cite{NeumeyerEtal2006} for linear models and \cite{AkritasKeilegom2001} for nonparametric regression. Under more restrictive assumptions one can further expand the remainder term as follows. Assume that $F$ is twice differentiable with density $F^\prime=f$ and bounded $f^\prime$ and further $|\hat\alpha-\alpha|+\|\hat\beta-\beta\|=o_{\mathbb{P}}(n^{-1/4})$, and $E\|X\|^2<\infty$. Then by Taylor's expansion
one obtains
$$\hat F_{\nt}(z)=F_{\nt}(z)+f(z)\left(\hat\alpha-\alpha+\langle E[X],\hat\beta-\beta\rangle\right)+o_{\mathbb{P}}(n^{-1/2}).$$ 
In models with intercept $\alpha$, where the estimator for $\alpha$ is chosen as $\hat\alpha=\overline{Y}_n-\langle \overline{X}_n,\hat\beta\rangle$  the remainder term is
$$R_n(z)=f(z)\left(\overline{\eps}_n+\langle \overline{X}_n-E[X],\hat\beta-\beta\rangle\right) +o_\mathbb{P}(n^{-1/2}).$$
(Here, $\overline{X}_n=n^{-1}\sum_{i=1}^n X_i$ and analogous for $\overline{Y}_n$ and $\overline{\eps}_n$.)
By Cauchy-Schwarz-inequality and the central limit theorem for $n^{1/2}(\overline{X}_n-E[X])$ one obtains that the dominating part of the remainder term is $f(z)\overline{\eps}_n$. This is the same as in 
homoscedastic finite-dimensional linear models with intercept and  nonparametric regression models. 
Note that $\alpha$ and $\beta$ are identifiable if the kernel of the covariance operator of the covariate $X$ is $\{0\}$. But often functional linear models without intercept are considered in the literature. So in our model assume $\alpha=\hat\alpha=0$. Then the remainder term is
$$R_n(z)=f(z)\langle E[X],\hat\beta-\beta\rangle +o_{\mathbb{P}}(n^{-1/2}),$$
and $\langle x,\hat\beta-\beta\rangle$ for fixed $x\in \mathcal{H}$ typically has a slower rate than $n^{-1/2}$, see \cite{CardotEtal2007}, \cite{ShangCheng2015}, \cite{YeonEtal2023}. If one assumes $E[X]=0$, then this problematic term cancels (similar as e.g.\ for centered ARMA-processes, see \cite{Bai1994}), but otherwise $f(z)\langle E[X],\hat\beta-\beta\rangle$ will dominate the asymptotic distribution of the process $(\hat F_{\nt}(z)-F(z))_{t\in [0,1],z\in\R}$. For our change-point test this dominating term vanishes.  The same holds  when estimating the conditional copula of the response in multidimensional functional linear models, given the covariate, see Theorem 5 in \cite{NeumeyerOmelka2024}. But e.g.\ for goodness-of-fit testing the remainder term would be of relevance. We consider goodness-of-fit tests for the error distribution in the different cases explained above in future work. 
With the derived expansion for residual empirical distribution functions one can also develop other tests for the error distribution as  e.g.\ for symmetry, or equality of error distributions in different models, see e.g.\
\cite{PierceKopecky1979}, 
\cite{NeumeyerEtal2005}, 
\cite{PardoFernandez2007}, among many others, in the cases of regression models with finite-dimensional covariates.

\begin{appendix}

\section{Proofs}

For ease of notation let $(X,Y,\eps)$ be some generic random variable with the same distribution as $(X_1,Y_1,\eps_1)$ under the null, but independent from the sample $(X_i,Y_i)$, $i=1,\dots,n$. Let $P$ denote the distribution of $(X,\eps)$. Further let $E_X$ denote the expectation with respect to $X$, which in the context below is the conditional expectation given $(X_i,Y_i)$, $i=1,\dots,n$.

  The proofs of Theorem \ref{theo1} and Lemma \ref{lem1} are similar as a part of the proof of Theorem 5 in \cite{NeumeyerOmelka2024}, but under less restrictive assumptions.

\subsection{Proof of Theorem \ref{theo1}}

From the Donsker-property in assumption (a.3) and Corollary 9.31 in \cite{Kosorok} it follows that
$$\left\{(x,e)\mapsto I\{e\leq z+a+\langle x,b\rangle\}-I\{e\leq z\}\mid z\in\R,a\in\R, b\in \mathcal{B}\right\}$$
is also $P$-Donsker. From Theorem 2.12.1 in \cite{VaartWellner} it follows that also the centered sequential process
$$H_n(t,z,a,b)=\frac{1}{\sqrt{n}}\sum_{i=1}^{\nt} \big( I\{\eps_i\leq z+a+\langle X_i,b\rangle\}
-I\{\eps_i\leq z\}-E[F(z+a+\langle X,b\rangle]+F(z)\big),
$$
indexed in $t\in[0,1]$, $z\in\mathbb{R}$, $a\in\mathbb{R}$, $b\in\mathcal B$,
converges weakly to a centered Gaussian process. 
Thus the process $H_n$ is asymptotically stochastic equicontinuous with respect to the semi-metric 
$$\rho((t_1,z_1,a_1,b_1),(t_2,z_2,a_2,b_2))=|t_1-t_2|+\Var( I\{\eps\leq z_1+a_1+\langle X,b_1\rangle\}-I\{\eps\leq z_2+a_2+\langle X,b_2\rangle\}),$$
see \cite{VaartWellner}, problem 2, p.\ 93, and section 2.12.  
In particular we need 
\begin{eqnarray*}
\rho((t,z,a,b),(t,z,0,0))
&=& \,\Var\left(I\{\eps\leq z+a+\langle X,b\rangle\}
-I\{\eps\leq z\}\right)\\
&\leq & \left|E\left[F(z+a+\langle X,b\rangle)-F(z)\right]\right|\\
&\leq& cE\left[|a+\langle X,b\rangle|^\gamma\right]\;\leq\; c\big(|a|+\|b\|E\|X\|\big)^\gamma
\end{eqnarray*}
by Cauchy-Schwarz and Jensen's inequality. 
Now setting $a=\hat\alpha-\alpha$ and $b=\hat\beta-\beta$ we obtain convergence to zero in probability by assumption (a.1). Thus from asymptotic stochastic equicontinuity of the process $H_n$, and $H_n(t,z,0,0)=0$ we obtain that 
$$\sup_{t\in[0,1],z\in \R}|H_n(t,z,\hat\alpha-\alpha,\hat\beta-\beta)|=o_{\mathbb{P}}(1),$$
which means that
\begin{eqnarray}\label{Entwicklung-hatF_nt}
\frac{\nt}{\sqrt{n}}\hat F_{\nt}(z)&=&\frac{\nt}{\sqrt{n}}F_{\nt}(z)+\frac{\nt}{\sqrt{n}}\left(E_X[F(z+\hat\alpha-\alpha+\langle X,\hat\beta-\beta\rangle)]-F(z)\right)\\
&&{}+o_{\mathbb{P}}(1)\nonumber
\end{eqnarray} 
uniformly in $t\in[0,1],z\in\R$, where
$F_{\nt}$ was defined in (\ref{F_nt}) and is based on the true errors. In particular for $t=1$ we have
$$\hat F_{n}(z)=F_{n}(z)+\left(E_X[F(z+\hat\alpha-\alpha+\langle X,\hat\beta-\beta\rangle)]-F(z)\right)+o_{\mathbb{P}}(n^{-1/2})$$ 
uniformly in $z\in\R$.
From those expansions we obtain 
\begin{eqnarray}\nonumber
\hat G_n(t,z)&=&\frac{\nt}{\sqrt{n}}\hat F_{\nt}(z)-
\frac{\nt}{\sqrt{n}}\hat F_n(z)\\
\nonumber
&=&\frac{\nt}{\sqrt{n}} F_{\nt}(z)-
\frac{\nt}{\sqrt{n}} F_n(z)+o_{\mathbb{P}}(1)\\
&=& G_n(t,z)+o_{\mathbb{P}}(1)\label{hat G_n=G_n}
\end{eqnarray}
uniformly in $t\in[0,1],z\in\R$.
\hfill $\Box$

\subsection{Proof of Lemma \ref{lem1}}

Let $\epsilon>0$ and let 
\[
 \big[b_i^L,b_i^U \big], \ i=1,\dotsc,N(\epsilon)=O\big(\exp(\epsilon^{-2/(k\gamma)})\big)
\]
be brackets for $\mathcal{B}$ of $\| \cdot\|$-length $\epsilon^{2/\gamma}$ (see assumption (a.3)). Now for  $b\in [b_i^L,b_i^U]$ the indicator function $I\{e\leq v + \langle x,b\rangle  \}$ is contained in the bracket 
\[
\big[\,I\{e\leq v+\langle xI\{x\geq 0\},b_i^L\rangle+\langle xI\{x<0 \},b_i^U\rangle  \},\, 
I\{e\leq v+\langle xI\{x\geq 0\},b_i^U\rangle+\langle xI\{x< 0\},b_i^L\rangle  \}
\big]
\] 
for each $v \in\mathbb{R}$. Further the above  
bracket has $L^2(P)$-length
\begin{eqnarray*}
&& \Big(E  \big[I\{\eps \leq v+\langle XI\{X\geq 0\},b_i^U\rangle+\langle XI\{X< 0\},b_i^L\rangle\}\\
&&\qquad{}
-I\{\eps \leq v+\langle XI\{X\geq 0\},b_i^L\rangle+\langle XI\{X< 0\},b_i^U\rangle\}\big]^2\Big)^{1/2}\\
&\leq& \left( E\big[F\big(v+\langle X,b_i^U\rangle \big) - 
 F\big(v+\langle X,b_i^L\rangle \big)\big] \right)^{1/2}\\
&\leq& \left( E\big[c|\langle X,b_i^U-b_i^L\rangle|^\gamma \big] \right)^{1/2}\\
&\leq &  c^{1/2}(E\|X\|)^{\gamma/2}\|b_i^U-b_i^L\|^{\gamma/2} \\
& =&O(\epsilon),
\end{eqnarray*}
by assumption (a.2), Cauchy-Schwarz and Jensen's inequality. 
 Similar to the proof of Lemma 1 in 
    \cite{AkritasKeilegom2001}
one obtains an upper bound $O(\epsilon^{-2}\exp(\epsilon^{-2/k}))$ for the  $L^2(P)$-bracketing number of the class $\mathcal{F}$. Thus $\mathcal{F}$ is $P$-Donsker by the bracketing integral condition in Theorem~19.5 of \cite{Vaart}. \hfill $\Box$

\subsection{Proof of Lemma \ref{lem-consistency}}

 To show the assertion for $\hat F_{k_n^*}$ we use the arguments as in the proof of Theorem \ref{theo1} for the process $H_n$, but based on the iid sample $(X_1,Y_1),\dots,(X_{k_n^*},Y_{k_n^*})$ before the change. Then as in the proof of Theorem \ref{theo1} asymptotic stochastic equicontinuity of the process $H_{k_n^*}$ holds and thus
 $$\sup_{t\in[0,1],z\in\mathbb{R}}|H_{k_n^*}(t,z,\hat\alpha-\alpha,\hat\beta-\beta)|=o_{\mathbb{P}}(1).$$
 Here, $\hat\alpha$ and $\hat\beta$ depend on the whole sample and assumption (a.1) is used. Thus as in equation (\ref{Entwicklung-hatF_nt}) we obtain
 $$\hat F_{k_n^*}(z)=F_{k_n^*}(z)+E_{X_1}[F_1(z+\hat\alpha-\alpha+\langle X_1,\hat\beta-\beta\rangle)]-F_1(z)+o_{\mathbb{P}}\big(\frac{1}{\sqrt{n}}\big)$$ 
 uniformly in $z\in\mathbb{R}$. By the classical Glivenko-Cantelli result $F_{k_n^*}$ converges uniformly almost surely to $F_1$. By assumptions (a.1), (a.2)' and $E\|X_1\|<\infty$ the remainder term is $o_{\mathbb{P}}(1)$ and  the assertion for $\hat F_{k_n^*}$ follows.  
 
 The assertion for $\tilde F_{k_n^*}$ can be shown analogously.
\hfill $\Box$

\subsection{Proof of Lemma \ref{lem-est-consistency}}

First note that
$$\hat\vartheta_n\in \arg\max_{t\in[0,1]}\Big(\sup_{z\in\mathbb{R}}|\hat G_n(t,z)|\Big)=\arg\max_{t\in[0,1]}\bigg(\sup_{z\in\mathbb{R}}\left|\frac{\hat G_n(t,z)}{n^{1/2}} \right|\bigg).$$
Further it holds
\begin{align*}
\frac{\hat G_n(t,z)}{n^{1/2}}  =&\ \frac{\nt \left(n-\nt\right)}{n^{2}}\Big(\frac{1}{\nt}\sum_{i=1}^{\nt} I\{\hat\eps_i\leq z\}-\frac{1}{n-\nt}\sum_{i=\nt+1}^n I\{\hat\eps_i\leq z\}\Big)\\
=&\ \frac{\nt (n-\nt)}{n^{2}}\bigg(\frac{1}{\nt}\sum_{i=1}^{\nt\wedge \lfloor n\vartheta^*\rfloor} I\{\hat\eps_i\leq z\}+I\{t> \vartheta^*\}\frac{1}{\nt}\sum_{i=\lfloor n\vartheta^*\rfloor+1}^{\nt} I\{\hat\eps_i\leq z\}\\
&\qquad-\frac{1}{n-\nt}\sum_{i=\nt\vee\lfloor n\vartheta^*\rfloor+1}^n I\{\hat\eps_i\leq z\}-I\{t< \vartheta^*\}\frac{1}{n-\nt}\sum_{i=\nt+1}^{\lfloor n\vartheta^*\rfloor} I\{\hat\eps_i\leq z\}\bigg)\\
=& \ \frac{\nt \left(n-\nt\right)}{n^{2}}\bigg(\frac{\nt\wedge \lfloor n\vartheta^*\rfloor}{\nt}F_1(z)+I\{t> \vartheta^*\}\frac{\nt- \lfloor n\vartheta^*\rfloor}{\nt}F_2(z)\\
&\qquad-\frac{n-\nt\vee \lfloor n\vartheta^*\rfloor}{n-\nt}F_2(z)-I\{t< \vartheta^*\}\frac{\lfloor n\vartheta^*\rfloor-\nt}{n-\nt} F_1(z)\bigg)+o_{\mathbb{P}}(1),
\end{align*}
since  we have
\begin{align}\nonumber
&\ \sup_{t\in[0,\vartheta^*]}\sup_{z\in\er}\frac{\nt}n\bigg|\frac{1}{\nt}\sum_{i=1}^{\nt} I\{\hat\eps_i\leq z\}-F_1(z)\bigg|\\
\leq&\ \sup_{t\in[0,\vartheta^*]}\sup_{z\in\er}\underbrace{\frac{\nt}{\lfloor n\vartheta^*\rfloor}\bigg|\frac{1}{\nt}\sum_{i=1}^{\nt} I\{\hat\eps_i\leq z\}-\frac{1}{\lfloor n\vartheta^*\rfloor}\sum_{i=1}^{\lfloor n\vartheta^*\rfloor} I\{\hat\eps_i\leq z\}\bigg|}_{=\frac 1{\lfloor n\vartheta^*\rfloor^{1/2}}\tilde G_{\lfloor n\vartheta^*\rfloor}(t,z)}\label{Term1}\\
&+\frac{\nt}{\lfloor n\vartheta^*\rfloor}\bigg|\frac{1}{\lfloor n\vartheta^*\rfloor}\sum_{i=1}^{\lfloor n\vartheta^*\rfloor} I\{\hat\eps_i\leq z\}-F_1(z)\bigg|\label{Term2}\\
=&\ o_{\mathbb{P}}(1). \nonumber
\end{align}
Here we have used Lemma \ref{lem-consistency} for the  term (\ref{Term2}). Further $\tilde G_{\lfloor n\vartheta^*\rfloor}$ is defined as $\hat G_n$ based on the iid-sample $(X_1,Y_1),\dots,(X_{k_n^*},Y_{k_n^*})$, but where the residuals are built with $\hat\alpha$, $\hat\beta$ based on the whole sample. With the same argument as in the proof of Theorem \ref{theo1} it holds that
$$\tilde G_{\lfloor n\vartheta^*\rfloor}(t,z)=G_{\lfloor n\vartheta^*\rfloor}(t,z)+o_{\mathbb{P}}(1)$$
uniformly in $t\in[0,1]$, $z\in\mathbb{R}$, see (\ref{hat G_n=G_n}), 
and thus the term (\ref{Term1}) is $o_{\mathbb{P}}(1)$. 

Analogously one can show that
$\sup_{t\in[\vartheta^*,1]}\sup_{z\in\er}\frac{n-\nt}n\big|\frac{1}{n-\nt}\sum_{i=\nt+1}^n I\{\hat\eps_i\leq z\}-F_2(z)\big|=o_{\mathbb{P}}(1)$.

Thus, it holds uniformly in $t\in[0,1]$
\begin{align*}
\frac{\hat G_n(t,z)}{n^{1/2}} 
=&\ I\{t> \vartheta^*\}\Big(\frac{\lfloor n\vartheta^*\rfloor(n-\nt)}{n^{2}} (F_1(z)-F_2(z)) \Big)\\
&+I\{t\leq \vartheta^*\} \Big( \frac{\nt(n-\lfloor n\vartheta^*\rfloor)}{n^{2}} (F_1(z)-F_2(z))\Big)+o_{\mathbb{P}}(1)\\
=&\ \Big(I\{t> \vartheta^*\}\vartheta^*(1-t)+I\{t\leq \vartheta^*\} t(1-\vartheta^*)\Big)(F_1(z)-F_2(z))+o_{\mathbb{P}}(1).
\end{align*}
The assertion then follows by Theorem 2.12 in \cite{Kosorok} as $\vartheta^*$ is well-separated maximum of $t\mapsto I\{t> \vartheta^*\}\vartheta^*(1-t)+I\{t\leq \vartheta^*\} t(1-\vartheta^*)$.

\hfill$\Box$

\end{appendix}

\bibliographystyle{apalike}
\bibliography{Bib}

\end{document}